\documentclass{article}
\usepackage{graphicx}

\newcommand{\pic}[4]{\vspace{1ex}\setlength{\unitlength}{1cm}
\begin{picture}(0,#3)(5,.5)
\put(#2){\includegraphics[#4]{#1.eps}}
\end{picture}\vspace{1ex}}

\newtheorem{proposition}{Proposition}[section]

\newtheorem{lemma}[proposition]{Lemma}
\def\proof{\textit{Proof. }}

\font\slashfont=msbm10 at 12pt
\font\slashfontsmall=msbm10 at 10pt
\def\C{\mbox{\slashfont C}}

\def\D{\mbox{\slashfont D}}
\def\Dsmall{\mbox{\slashfontsmall D}}
\def\DS{\displaystyle}

\def\R{\mbox{\slashfont R}}
\def\S{{\mathcal{S}}}
\def\qed{\hskip1em\raise3.5pt\hbox{\framebox[2mm]{\ }}}
\def\re{\,\!\mbox{Re}\,}
\def\im{\,\!\mbox{Im}\,}
\def\scq{s.c.q.}
\let\conj\overline
\def\str{\raisebox{1.5ex}{\rule{0pt}{1ex}}}
\newcommand{\strd}[1]{\left#1\str\right.}
\newcommand{\X}[1]{X^{(#1)}}
\newcommand{\Xt}[1]{\widetilde X^{(#1)}}

\begin{document}

\begin{center} 
{\LARGE
 Conformal Mapping\\ of Right Circular Quadrilaterals \\[1ex]
}
{\large Vladislav V. Kravchenko\footnote{Partially supported by CONACyT grant 50424}  
\\R. Michael  Porter\footnote { Partially supported by CONACyT grant 80503}
}

\medskip
Departamento de Matem\'aticas,\\
Cinvestav del I.P.N, Campus Quer\'etaro,\\ Apartado Postal 1-798, Arteaga \#5, Col.\ Centro,\\ Santiago de Quer\'etaro, Qro. 76001, Mexico

\medskip
\today
\end{center}

 \begin{minipage}{\textwidth}

\small 

\noindent Abstract. 
We study conformal mappings from the unit disk to circular-arc
quadrilaterals with four right angles.  The problem is reduced to a
Sturm-Liouville boundary value problem on a real interval, with a
nonlinear boundary condition, in which the coefficient functions
contain the accessory parameters $t$, $\lambda$ of the mapping
problem.  The parameter $\lambda$ is designed in such a way that for
fixed $t$, it plays the role of an eigenvalue of the Sturm-Liouville
problem.  Further, for each $t$ a particular solution (an elliptic
integral) is known a priori, as well as its corresponding spectral
parameter $\lambda$.  This leads to insight into the dependence of the
image quadrilateral on the parameters, and permits application of a
recently developed spectral parameter power series (SPPS) method for
numerical solution.  Rate of convergence, accuracy, and computational
complexity are presented for the resulting numerical procedure, which
in simplicity and efficiency compares favorably with previously known
methods for this type of problem.



\medskip
   
 \end{minipage}

\setcounter{section}{-1}
\section{Introduction}\label{sec:intro}

A \textit{symmetric circular quadrilateral} (\scq) is a Jordan curve
$P$ in $\C$ formed of four circular arcs (or straight segments) with
all four internal angles equal to $\pi/2$.  Will assume that the
vertices of $P$ are situated in positions of the form $\pm
A,\pm\conj{A}$, where $\re A>0$, $\im A>0$.

Let $D$ be a plane domain containing the origin and bounded by an
\scq\ The set of conformal types of such domains (considering the
vertices as distinguished points) forms a two-dimensional space in a
natural way.  We consider conformal mappings $f\colon\D\to{}D$ from
the unit disk $\D$ to $D$.  Since $f$ extends continuously to the
boundary $P$, there is a unique value $t\in[0,2\pi]$ such that
$f(e^{it})=A$; We will generally assume that $f(0)=0$ and $f'(0)>0$,
which implies $0<t<\pi/2$.

The more general problem of mapping the disk to circular-arc polygon
domains is treated in \cite{BG}, \cite[Chapter 4]{DT}, \cite[Chapter
  16]{Hen}, \cite[Section 17.6]{Hi}, \cite{Ho}, \cite[Chapter 5]{Neh}.
In particular, it is well known that the Schwarzian derivative of a
conformal mapping of $\D$ onto a circular-arc polygon is a rational
function of degree two. In this article we will develop further the
work of P.\ Brown on the accessory-parameter problem for \scq s.  We
follow much of the notation and copy several equations from \cite{Br1}
where it is verified that due to the symmetries of \scq s, the
Schwarzian derivative $\S_f$ of $f$ is of the specific form
(\ref{eq:Sforig}) below. This function is determined by two real
parameters $t,s$ and the relationship (\ref{eq:rho}) must hold.  As a
partial converse, it is well known that due to the intimate
relationship between the Schwarzian derivative and curvature, if the
Schwarzian derivative of a holomorphic function $f$ defined in $\D$ is
of the form (\ref{eq:Sforig}) and if (\ref{eq:rho}) holds, then $f$ is
a local homeomorphism onto a (not necessarily schlicht) domain bounded
by a (not necessarily simple) right circular-arc polygon.
 
A basic question is the following. Given $P=\partial D$ (for example
by specifying $A$ and also the radius or the midpoint of one of the
edges of $P$), to find the parameters $t,s$ of the Schwarzian
derivative of $f$. In \cite{Br1} Brown looked first at a simpler
question, fixing $t$, normalizing $f'(0)=0$ and calculating the
remaining parameter $s$ corresponding to a geometric characteristic of
$P$, such as the radius of one of its edges.  This reduces the problem
to one real dimension.  As with most methods which have been developed
for conformal mapping of circular-arc polygons, in \cite{Br1} the
Schwarzian differential equation is solved numerically for trial
values of $s$, the corresponding geometric characteristics of
$P=f(\partial\D)$ are calculated, and the process is repeated until a
sufficiently close $s$ is located.  Here we apply a technique
developed in \cite{KP} for dealing with Sturm-Liouville problems,
the spectral parameter power series (SPPS) method,
which permits a direct calculation in the sense that once certain
auxiliary parameters are calculated for given $t$, one may evaluate
the solution corresponding to any desired $s$ without resorting to
further integrations.

In a second article \cite{Br2} the full two-parameter problem is
addressed.  The problem is formulated with the normalization $f(1)=1$,
and the data is given in terms of the radius $1/\kappa_1$ of the right
edge of $P$ and the midpoint $p_2$ of the upper edge.  Brown's solution
involves a table of previously calculated values of these parameters
in terms of $(t,s)$.  To apply it one looks for values reasonably
close to the desired geometric parameters in this table, and then
applies an iterative process to approximate the sought-after $(t,s)$
to the desired accuracy.  In this paper we apply our solution of the
one-parameter problem, which is quite rapid, to this two-parameter
problem in an iterative way. In part due to properties of an
equivalent parameter $\lambda$ which we use in place of $s$, our
method does not require consultation of a table of prior values.

In the next two sections we set up a Sturm-Liouville boundary value
problem whose solution relates the accessory parameter $\lambda$ to
the curvature $\kappa$ of the right edge of $P$. It is seen that
$\lambda$ is a spectral parameter in a boundary value problem. In
Section \ref{sec:DegSol} we describe the so-called canonical mapping
$f_\infty$ of $\D$ to a rectangle, and identify its parameter
$\lambda_\infty$ and the corresponding eigenfunction $y_\infty$ of the
Sturm-Liouville problem, which are used in the application of the SPPS
method which is summarized in Section \ref{sec:II} and then
applied in Section \ref{sec:curvature}   to represent $\kappa$ as a power
series in $\lambda$ for fixed $t$.  An algorithm for the one-parameter
problem $\kappa\mapsto\lambda$ for fixed $t$ is presented in Section
\ref{sec:1par}, and for the two-parameter problem
$(\kappa_1,p_2)\mapsto(t,\lambda)$ in Section \ref{sec:2par}.

\section{Sturm-Liouville equation \label{sec:SLeq}}

We begin by setting up the classical second-order differential
equation which governs the conformal mapping to \scq s. The Schwarzian
derivative of a holomorphic $f$, namely
$S_f=(f''/f')'-(1/2)(f''/f')^2$, is again holomorphic when $f'$ does
not vanish.  For the conformal mapping of $\D$ to an \scq, formula
\cite[(3)]{Br1}, which we will refer to as (B3), says
\begin{eqnarray} \label{eq:Sforig}
 \S_f(z) &=& \frac{3}{8}\left(
   \frac{1}{(z-e^{it})^2} +
   \frac{1}{(z+e^{it})^2} +
   \frac{1}{(z-e^{-it})^2} +
   \frac{1}{(z+e^{-it})^2} \right)   \nonumber\\
   && +\ \ \left(
   \frac{c}{z-e^{it}} -
   \frac{c}{z+e^{it}} +
   \frac{c}{z-e^{-it}} -
   \frac{c}{z+e^{-it}} \right) 
\end{eqnarray}
in which the parameter $c\in\C$ is subject to one restraint as
follows.  Write $c=\rho e^{is}$. Then (B8)
\begin{equation} \label{eq:rho}
  \rho = \frac{-3}{8\cos(s+t)}
\end{equation}
where further we may assume (B9) $\pi/2-t<s<3\pi/2-t$.

We now introduce a new parameter equivalent to $s$.  By
(\ref{eq:rho}), the quantity $\rho>0$ satisfies in fact $\rho\ge3/8$.
Since $s=\arccos(-3/(8\rho))-t$ where
$\pi/2<\arccos(-3/(8\rho))<3\pi/2$, the parameter $c$ in
(\ref{eq:Sforig}) can be expressed as follows:
\begin{eqnarray} \label{eq:c}
 c&=& \rho e^{is} = \rho\left(
      \cos ((\arccos\frac{-3}{8\rho}) - t ) 
  + i \sin ((\arccos\frac{-3}{8\rho}) - t ) \right) \nonumber\\
 &=& -\frac{3}{8}(\cos t - i \sin t) +
       \varepsilon\sqrt{\rho^2-(3/8)^2} \,(\sin t +i\cos t) \nonumber \\
 &=& -\frac{3}{8} e^{-it} + \lambda i  e^{-it} 
\end{eqnarray}
where we write
\[   \lambda = \varepsilon \sqrt{\rho^2-(3/8)^2} 
\]
with
\[  \varepsilon= \left\{ \begin{array}{ll}
         -1,\quad&s<\pi-t,\\
         0,& s=\pi-t,\\
         1,& \pi-t<s. 
\end{array} \right.
\]
At the particular value $s=\pi-t$ we have $\rho=\infty$,
and also $\rho\to\infty$ as $s$ tends to either extreme value
$\pi/2-t$ or $3\pi/2-t$. Further, it is easily seen that
\begin{equation}
  \lambda =   \frac{3}{8}\tan(s+t).
\end{equation}

Everything we will do depends on the fact, manifested in (\ref{eq:c}),
that the parameter $c$ is a linear polynomial in $\lambda$.
The same then holds for $S_f(z)$. Indeed, abbreviating $a=e^{it}$, we
find that
\begin{eqnarray*} 
 \S_f(z) &=& \left( \frac{3}{8} 
\strd{(}
    \frac{1}{(z-a)^2} +
    \frac{1}{(z+a)^2} +
    \frac{1}{(z-\conj{a})^2} +
    \frac{1}{(z+\conj{a})^2} 
\strd{)} 
\right. \nonumber \\
  && \ \   \left. - \frac{3}{4} \strd{(}
\frac{1}{z^2-a^2} + \frac{1}{z^2-\conj{a}^2} \strd{)} \right) \nonumber\\
  && \  + 2\lambda i \strd{(}
     \frac{1}{z^2-a^2} - \frac{1}{z^2-\conj{a}^2} \strd{)} \nonumber \\
  &=& 2\psi_0(z) - 2 \lambda \psi_1(z)
\end{eqnarray*}
where  
\begin{eqnarray} \label{eq:psi}
  \psi_0(z) &=& \frac{3}{2}\left(
      \frac{a^2}{(z^2-a^2)^2} +
      \frac{\conj{a}^2}{(z^2-\conj{a}^2)^2} \right) \nonumber\\
  \psi_1(z) &=& -i \left(
      \frac{1}{z^2-a^2} -\frac{1}{z^2-\conj{a}^2} \right)
\end{eqnarray}

From now on we will regard $\S_f$ as parametrized by $t,\lambda$
instead of $t,s$. The functions $\psi_0,\psi_1$ depend also on $t$, so
when we require more clarity we will write $S_f(z)=R_{t,\lambda}(z)$
where
\begin{equation} \label{eq:R}
     R_{t,\lambda}=2\psi_0-2\lambda\psi_1.
\end{equation}
It is a textbook fact \cite{Hi,Neh} that all holomorphic functions having
Schwarzian derivative equal to $R_{t,\lambda}$ are quotients
$y_2/y_1$
where $y_1,y_2$ are linearly independent solutions of the ordinary
differential equation $2y''+R_{t,\lambda}\,y=0$ in $\D$, i.e., solutions
of
\begin{equation} \label{eq:SLeq}
 y'' +  \psi_0 y = \lambda \psi_1 y .
\end{equation}
It is our purpose to take advantage of the particularly nice form
of equation (\ref{eq:SLeq}) and its relation to spectral theory.

\section{Boundary Value Problem \label{sec:BVP}}

In \cite{Br1} Brown describes in detail the possible relations among
the edges of \scq s.  In particular, for a generic \scq, the two
circles containing one pair of opposite edges will intersect, while
the two circles containing the remaining two edges are disjoint.
There are two degenerate cases: when one pair of opposite edges lie on
the same circle, or when all four edges are straight segments (i.e.,
when $P$ is rectangular).  For the generic case, following \cite{Br1}
we will take the right and left edges to lie in disjoint circles of
common radius $|r|$, where $r$ is given (B14) by
\begin{equation} \label{eq:radius}
 r = \frac{f'(1)^2}{f'(1)+f''(1)} .
\end{equation}
We will normalize $f$ by 
\begin{equation} \label{eq:f0}
 f(0)=0,\quad f'(0)=1,\quad f''(0)=0
\end{equation}
recalling that the Schwarzian equation has three degrees of freedom.
(This means that the image may need to be rescaled, i.e.,
$P=\mbox{constant}\cdot f(\partial\D)$).

Let $y$ be a solution of (\ref{eq:SLeq}) in $\D$ normalized by
$y(0)=1$, $y'(0)=0$. Then the normalized mapping $f$ is equal to the
indefinite integral
\begin{equation} \label{eq:fintegral}
 f(z) = \int_0^z \frac{1}{y(\zeta)^2} \,d\zeta.
\end{equation}
Indeed, (\ref{eq:fintegral}) yields $f'=y^{-2}$ and $f''=-2y^{-3}y'$, 
so $f$ satisfies the required normalization $(\ref{eq:f0})$.
Now from (\ref{eq:radius}) we have
\[ 0 = r\strd{(}f'(1)+f''(1)\strd{)}-f'(1)^2 
     = r\strd{(}y(1)^{-2} - 2y(1)^{-3}y'(1)\strd{)} - y(1)^{-4}
\]
so the curvature of the right edge of the image is given by
\begin{equation} \label{eq:kappa}
 \kappa = \frac{1}{r} =
   y(1)^2 - 2y(1)y'(1).
\end{equation}
Consequently the boundary value problem to solve is for $y(z)$ with
real $z\in[0,1]$,
\begin{equation} \label{eq:bvp}
 \left\{ \ \ \parbox{.5\textwidth}{
  $\DS y'' + \psi_0 y = \lambda \psi_1 y$,\\[1ex]
  $\DS y(0) = 1$, \\[1ex]
  $\DS y'(0) = 0$, \\[1ex]
  $\DS y(1)^2 - 2y(1)y'(1) = \kappa$. 
	} \right.
\end{equation}
The data for this problem are $t$ and $\kappa$, and one seeks
$\lambda$ for which there is a solution to the system (\ref{eq:bvp}),
where we recall that the dependence on $t$ is through $\psi_0$ and
$\psi_1$.  Among such $\lambda$ one of course wants the value for
which the corresponding mapping $f$ is injective. In \cite{Br1} there
are many details regarding the values of $t,s$ for which this holds,
and one can apply these facts to $\lambda$.

The nonlinear nature of the boundary condition at $z=1$ in
(\ref{eq:bvp}) will be dealt with in Section \ref{sec:curvature}.

\section{The Degenerate S.C.Q.s \label{sec:DegSol}}

Again following \cite{Br1}, 
let $O(s)\in\R$ denote the center of the circle containing the
rightmost edge of the \scq\ $f(\D)$.  In \cite{Br1} it is shown that
$O(s)$ is a monotone function of $s$ for each fixed $t$.  Two special
values $s_0$, $s_\infty$ are singled out, with $O(s_0)=0$ (called the
\textit{root mapping}) and $O(s_\infty)=\infty$ (the \textit{canonical
mapping}), for which $f(\D)$ has the degenerations mentioned at the
begining of Section \ref{sec:BVP}.  The root mapping and its parameter
$s_0$ are worked out in detail.  We will use $s_\infty$ here because
it is simpler; in fact, the canonical mapping $f_\infty$ is an
elliptic integral of the first kind:
\begin{equation} \label{eq:finfty}
  f_\infty(z) = \int_0^z (1-2(\cos 2t)\zeta^2 + \zeta^4)^{-1/2}\,d\zeta.
\end{equation}
Here $f_\infty(\D)$ is a rectangle, $f_\infty([-1,1])\subseteq\R$,
and the values $k_1=f_\infty(1)$ and  $k_2=-if_\infty(i)$ are known
clasically as complete elliptic integrals.

Note that $f_\infty'$ never vanishes, and  $f_\infty'(0)=1$.
By direct calculation one finds that
\begin{eqnarray*}
 \S_{f_\infty} &=&  \left(\frac{f_\infty''}{ f_\infty'}\right)'
    -\frac{1}{2} \left(\frac{f_\infty''}{ f_\infty'}\right)^2  \\[1ex]
   &=&  \frac{2(\cos2t+(\cos^22t-3)z^2+(\cos2t)z^4)}{
              (1-2(\cos2t)z^2+z^4)^2} .
\end{eqnarray*}
However, by (B23),
\begin{equation} \label{eq:B23}
  \S_{f_\infty} = 
    \frac{2(\cos2t+(\cos^22t-3)z^2+(\cos2t)z^4)}{
              (1-2(\cos2t)z^2+z^4)^2} +
    \frac{\DS 3\frac{\cos(s-t)}{\cos(s+t)} - \cos2t}{
                 2(1-2(\cos2t)z^2+z^4)} .
\end{equation}
By comparison, the second term on the right side of (\ref{eq:B23})
must vanish,
\[   3\frac{\cos(s-t)}{\cos(s+t)} = \cos2t .
\]
This can be written many ways, for example
\begin{equation} \label{eq:sfromt}
 \tan s = (\cot t) \frac{\cos2t-3}{\cos2t+3}  
      = -(\cot t) \frac{ \frac{1}{\sin t}+\sin t}
                       { \frac{1}{\cos t}+\cos t}
\end{equation}
etc.  Given $t$, (\ref{eq:sfromt}) provides $s=s_\infty$ for the
canonical mapping.  With this in hand, we obtain the parameters
\begin{eqnarray*}
   \rho_\infty &=& \frac{-3}{8\cos(s_\infty+t)},\\
   \lambda_\infty &=& \varepsilon_\infty
                      \sqrt{\rho_\infty^2-\left(\frac{3}{2}\right)^2}\\
     &=& \frac{-3}{8}\tan(s_\infty+t) .
\end{eqnarray*}
With the aid of (\ref{eq:sfromt}) one can simplify this last equation to
\begin{equation}
   \lambda_\infty = \frac{1}{4}\cot 2t
\end{equation}

The Sturm-Liouville system (\ref{eq:bvp}) corresponding to the
canonical mapping with Schwarzian derivative $R_{t,\lambda_\infty}$
has the particular solution $y_\infty=(f_\infty)^{-1/2}$, specifically
\begin{equation} \label{eq:yinf}
  y_\infty(z) = (1 -2(\cos2t)z^2+z^4)^{\frac{1}{4}}
\end{equation}
Clearly this satisfies the normalizations $y_\infty(0)=1$ and
$y_\infty'(0)=0$, and for reference we note
\begin{equation}  \label{eq:yinfvals}
  y_\infty(1)=\sqrt{2\sin t},\quad\quad
  y_\infty'(1)=\sqrt{\frac{\sin t}{2}}.
\end{equation}
Since $f_\infty(\D)$ is a rectangle, the curvature of its right edge is
$\kappa_\infty=0$.  Further, since $f_\infty'$ never vanishes, it
follows that the function $y_\infty$ also never vanishes in $\D$; this
can also be seen from (\ref{eq:yinf}).

\section{Solution by Iterated Integrals \label{sec:II}}

We now apply the spectral parameter power series (SPPS) method for
solution of Sturm-Liouville problems, developed in \cite{KP}. This will
provide power series in $\lambda$ which represent the geometric
parameters of the \scq\ mapping problem.

Let two functions $q_0$, $q_1$ be given on $[0,1]$.  (All of
the following is valid for $q_n(z)$ defined in, say, $|z|<1$, but the
numerical method presented later will not involve complex values of
$z$.) The sequence $I_n$ of \textit{iterated integrals} determined by
the generating pair $(q_0,q_1)$ is defined as follows.  Let $I_1=1$
identically on $[0,1]$, and then recursively for $n\ge2$,
\begin{equation}
  I_n(z) = \int_0^z I_{n-1}(\zeta)\,q_{n-1}(\zeta)\,d\zeta
\end{equation}
where the indices in the generating pair are understood mod 2; i.e.,
$q_{n+2j}=q_n$ for $j=1,2,\dots$

The relationship between the iterated integrals and the basic equation
(\ref{eq:SLeq}) is given by the following result.

\begin{proposition}{\rm \cite{KP}} \label{prop:powersX}
 Let\/ $\psi_0$ and $\psi_1$ be given, and suppose that\/ $y_\infty$
is a function which does not vanish and which satisfies the 
ordinary differential equation
\[ y_\infty''+\psi_0\,y_\infty=\lambda_\infty\psi_1\,y_\infty
\] on $[0,1]$. Choose\/ $q_0=1/y_\infty^2$, $q_1=\psi_1\,y_\infty^2$
and define\/ $\X{n}$, $\Xt{n}$ to be the iterated integrals
determined by\/ $(q_0,q_1)$ and by\/ $(q_1,q_0)$ respectively.  Then for
every\/ $\lambda\in\C$ the functions
\begin{eqnarray*}
  y_1 = y_\infty \sum_{k=0}^\infty(\lambda-\lambda_\infty)^k
            \Xt{2k},\\ 
  y_2 = y_\infty \sum_{k=0}^\infty(\lambda-\lambda_\infty)^k
              \X{2k+1}\\
\end{eqnarray*}
are linearly independent solutions of the equation
\[ y''+\psi_0y=\lambda\psi_1y 
\]
on the interval\/ $[0,1]$.  Further, the series for $y_1$ and $y_2$ 
converge uniformly on $[0,1]$ for every $\lambda$.
\end{proposition}

\section{Variation of Curvature \label{sec:curvature}}

\begin{figure}[!b]
\pic{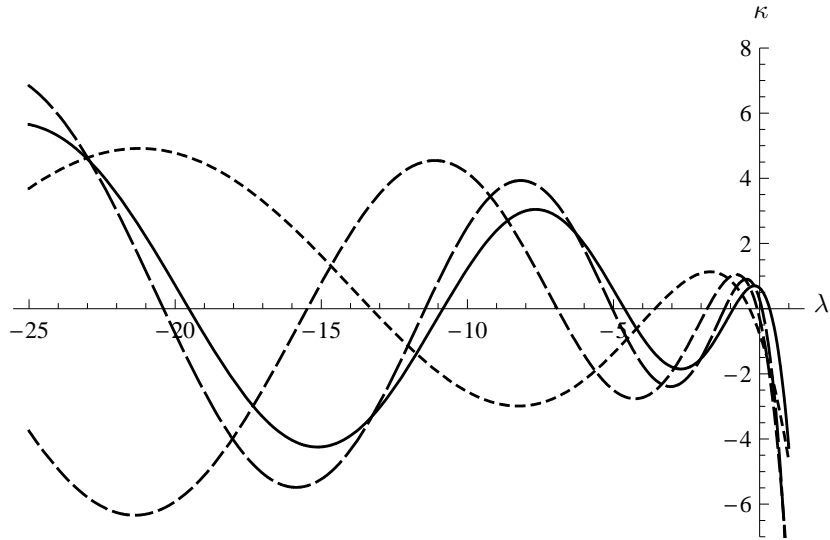}{5.5,0}{6}{}
\caption{\small Graph of $\kappa$ as a function of $\lambda$ for
  $t=0.1\pi, 0.2\pi,0.3\pi,0.4\pi$.
\label{fig:kapgraph} }
\end{figure}

To apply Proposition \ref{prop:powersX} to the boundary value problem
(\ref{eq:bvp}), we take $\psi_0,\psi_1$ as in (\ref{eq:psi}) and
$y_\infty$ as in (\ref{eq:yinf}). Note that
\[ y_1' = y_\infty'\sum_{k=0}^\infty(\lambda-\lambda_\infty)^k\Xt{2k} +
  \frac{1}{y_\infty}\sum_{k=1}^\infty (\lambda-\lambda_\infty)^k\Xt{2k-1}
\] 
with a similar formula for $y_2'$, so the initial values of the
particular solutions $y_1$, $y_2$ are given by
\begin{eqnarray*}
  y_1(0)=y_\infty(0)=1,  && y_1'(0)=y_\infty'(0)=0,\\
   y_2(0)=0,             &&
   y_2'(0)=\frac{1}{y_\infty(0)}=1 .
\end{eqnarray*}
This says that $y_1$ as given by Proposition \ref{prop:powersX} is
precisely the solution of the second order linear differential
equation which satisfies the two initial conditions of (\ref{eq:bvp})
at the point $z=0$.  Consequently problem (\ref{eq:bvp}) reduces to
satisfying the final boundary condition,
\begin{equation} \label{eq:curvecond}
  y_1(1)(y_1(1)-2y_1'(1)) = \kappa . 
\end{equation}
The numerical evaluation of the left side of (\ref{eq:curvecond}) is
readily accessible because it is a power series $\sum
a_n(\lambda-\lambda_\infty)^n$ in $\lambda$ whose coefficients $a_n$
are represented in terms of those of the series $y_1(1)=\sum
b_n(\lambda-\lambda_\infty)^n$ provided by Proposition
\ref{prop:powersX}.

In Figure \ref{fig:kapgraph} we show a graph of the power series 
(\ref{eq:curvecond}) for several values of $t$.  Numerical details
will be given later; here we wish to stress that once having
calculated the coefficients $a_n$ for a given $t$, one needs no longer
to solve a (Schwarzian or Sturm-Liouville) differential equation in
order to calculate $\kappa$ as a function of $\lambda$ for that value
of $t$.

Perhaps surprisingly, the calculation of (\ref{eq:curvecond}) can be
simplified much further.  Note that the $k$-th Taylor coefficient  of
$y_1(1)-2y_1'(1)$ is
\begin{eqnarray*}
 b_k - 2b'_k &=& y_\infty(1)\Xt{2k}(1)-2\left(
            y_\infty'(1)\Xt{2k}(1) +\frac{1}{y_\infty(1)}\Xt{2k-1}(1)
        \right) \\
   &=& \Xt{2k}(1)(y_\infty(1)-2y_\infty'(1))
         - \frac{2}{y_\infty(1)}\Xt{2k-1}(1).
\end{eqnarray*}
It follows from (\ref{eq:yinf}) that $y_\infty(1)-2y_\infty'(1)=0$,
which leaves
\[
   b_k - 2b'_k =  - \frac{2}{y_\infty(1)}\Xt{2k-1}(1).
\]
Thus 
\begin{equation} \label{eq:an}
  a_n = -2\sum_{k=0}^n \Xt{2k}(1) \, \Xt{2(n-k)-1}(1).
\end{equation}

From (\ref{eq:psi}),(\ref{eq:yinf}) it is also clear that
$y_\infty(z)\ge0$ and $\psi_1(z)\ge0$ for $z\in[0,1]$, and
consequently $\Xt{k}(z)\ge0$.  Thus we see that in fact all the Taylor
coefficients $a_n$ of $\kappa$ as a function of $\lambda$ are
negative, with the exception of $a_0=0$. 

\begin{figure}[!b]
\pic{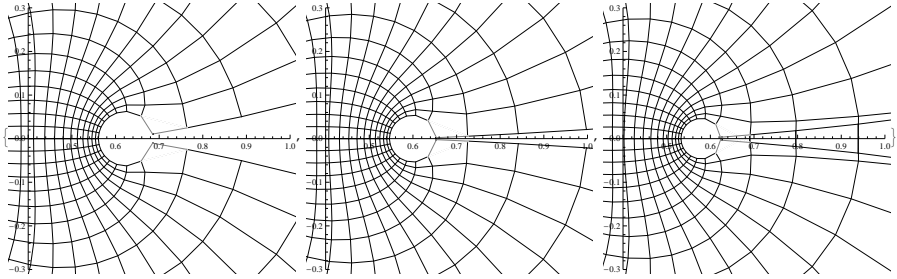}{5.2,.3}{4}{scale=.45}
\caption{\small Images $f(\Dsmall)$ for $t=\pi/4$ and
$\lambda=1.3,\ 1.4,\ 1.5$. The images are depicted near the point
$w=0.5$, where it is perceived that the rightmost figure is a
non-schlicht region. 
\label{fig:biglambda} }
\end{figure}

Note that $\kappa$ is never zero for $\lambda>\lambda_\infty$; in
fact, $d\kappa/d\lambda$ is strictly negative, and $\kappa\to-\infty$
as $\lambda\to\infty$.  Thus the radius of the rightmost edge of
$f(\D)$ tends to zero, as illustrated in Figure \ref{fig:biglambda}.  As
$\lambda\to\infty$, the arc $e^{i\theta}:-t<\theta<t$ of $f(\partial
D)$ is mapped to a curve tracing ever-greater numbers of full turns
around ever-smaller circles.

\begin{figure}[!t]
\pic{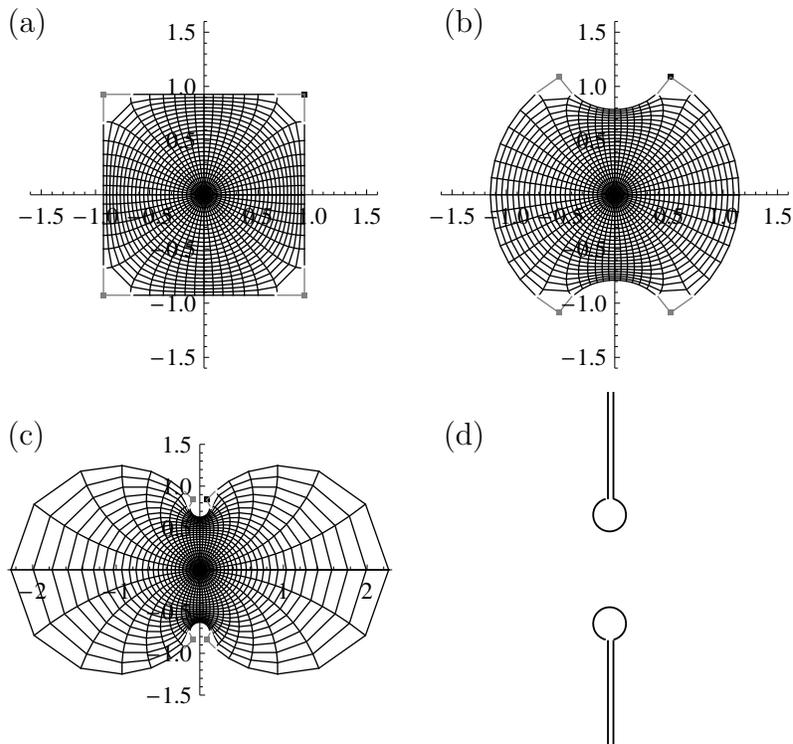}{6.,0}{9.2}{scale=1}
\caption{\small Images $f(\Dsmall)$ for $t=\pi/4$ and (a)
  $\lambda=0=\lambda_\infty=0$, (b) $\lambda=-0.32219$,
  $\lambda=-0.91570$ and (d) $\lambda=-1.43554$ (schematic drawing of
  boundary).  The value in (d) corresponds to second largest zero of
  $\kappa(\lambda)$, after $\lambda_\infty$.  The
  values for (b), (c) were chosen to give equal values $\kappa=0.8$.
 \label{fig:neglambda} }
\end{figure}

In Figure \ref{fig:neglambda} we bend the edge in the opposite
direction, taking $\lambda<\lambda_\infty$.  The largest zero of
$\kappa(\lambda)$ is $\lambda_\infty$ corresponding to the canonical
mapping.  The next zero gives a domain, exterior to the non-simple
curve depicted in Figure \ref{fig:neglambda}(d), covering the exterior
in the Riemann sphere of two circles joined by a straight segment.
Further zeros correspond to non-schlicht domains.

\section{Algorithm for the One-Parameter Problem \label{sec:1par}}

Here we solve the nonlinear problem (\ref{eq:bvp}) mentioned at the
end of Section \ref{sec:BVP}. Fix $0<t<\pi/2$.  Given $\kappa$, we
seek a parameter $\lambda$ so that the mapping $f$ with Schwarzian
derivative $\S_f=R_{t,\lambda}$, normalized by $f(0)=0$, $f'(0)=1$
will send $\partial\D$ to an \scq\ whose right edge has curvature
$\kappa$.  (If $\kappa<0$, we want the center of the edge arc to lie
to the left, i.e., in the ray $(-\infty,1)$.) We stress again that for
fixed $t$, once certain parameters have been determined, the relation
$\kappa\mapsto\lambda$ is of rapid calculation.

The algorithm given below summarizes results of the previous sections,
in particular formulas (\ref{eq:psi}), (\ref{eq:yinf}), and
(\ref{eq:an}). In the following, ``evaluate'' means to calculate at
$M+1$ evenly spaced values of $z$ from 0 to 1.  We will use $N+1$
terms in the power series in $\lambda$.  Observe that $\psi_0$ and
$\X{n}$ are not needed.

\medskip
\noindent\textit{One-parameter algorithm.}
\begin{enumerate}

\item Evaluate $y_\infty(z) = (1-2(\cos2t)z^2+z^4)^{1/4}$.

\item   Evaluate
\begin{eqnarray*} 
  \psi_1(z) &=& \frac{2\sin 2t}{z^4-2(\cos 2 t)z^2 + 1} .
\end{eqnarray*}

\item Evaluate the iterated integrals $\Xt{n}$ determined by
$(\psi_1y_\infty^2,\ 1/y_\infty^2)$ for $n=0,1,\dots,2N$.


\item
 Calculate the Cauchy product terms 
\[  a_n = -2\sum_{k=0}^n \Xt{2k}(1) \, \Xt{2(n-k)-1}(1).
\]
for $n=0,1,\dots,N$. Calculate also $\lambda_\infty=(1/4)\cot2t$.

\item For any desired $\kappa$, solve the polynomial equation
\[   \sum_{n=0}^{N}a_n(\lambda-\lambda_\infty)^n=\kappa
\]
for $\lambda-\lambda_\infty$, and then add $\lambda_\infty$ to
the result obtain $\lambda$.
\end{enumerate}

\medskip\noindent Then the mapping $f_{t,\lambda}$ will produce curvature
approximately $\kappa$ on the right and left edges.

\section{Two-Parameter Problem \label{sec:2par}}

We now turn to the two-parameter mapping problem addressed in
\cite{Br2}.  Given an \scq\ $P$ proposed to be the image of $\D$ under
$w=f(z)$, one does not know the value of $f'(0)$ directly from the
geometry of $P$, so instead one assumes that the right edge of $P$
passes through $w=1$.  Consider the mapping $g(z)= f(z)/f(1)$
which satisfies
\[ g(0)=0,\quad g(1)=1. 
\]
We let $\kappa_1$ and $\kappa_2$ denote the curvatures of the right
and upper edges of $P=g(\partial\D)$, with midpoints $p_1=1$ and
$p_2\in i\R^+$.  The basic mapping problem will be to find $(t,\lambda)$
given $(\kappa_1,p_2)$.

\medskip
\noindent\textit{Geometric parameters of $P$.}
To begin with, we take into account the values which we could readily
calculate if we were already given $(t,\lambda)$.  Write $w_1=f(1)$,
$w_2=f(i)$.  We have already seen how to calculate
$w_1=y_2(1)/y_1(1)$.  For $w_2$, observe that by
(\ref{eq:psi}),(\ref{eq:R}),
\begin{equation} \label{eq:Rsymm}
   R_{\pi/2-t,-\lambda}(\pm iz) = -R_{t,\lambda}(z).
\end{equation}
Further,  $w_2/i=f^*(1)$ where $f^*(z)=-if(iz)$, and Chain Rule for
Schwarzian derivatives gives us
\[    S_{f^*}(z) = -S_f(iz) = R_{t',\lambda'}(z).
\]
Thus $w_2$ can be obtained simply by using $(\pi/2-t,-\lambda)$ in
place of $(t,\lambda)$ in the calculation of $w_1$. This provides us
in turn with the desired values
\[   p_2 = \frac{w_2}{w_1}
\]
and   
\[  \kappa_1 = \kappa w_1
\]
where $\kappa$ denotes as previously the curvature of the right side
of $(1/w_1)P=f(\D)$.  As we have seen, for fixed $t$ these quantities are
power series in $\lambda$.

\medskip
\begin{figure}[!b]
\pic{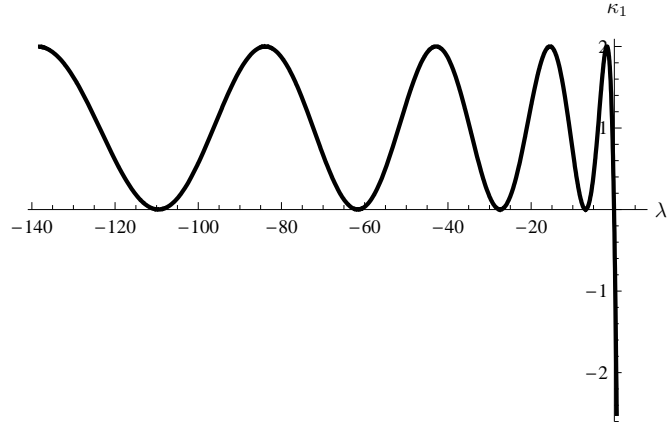}{6.5,.3}{4.9}{scale=.8}
\caption{\small Normalized curvature $\kappa_1$ as a function of
  $\lambda$ for $t=0.3\pi$.  As $\lambda$ passes through the local
  maxima at height $\kappa_1=2$, the boundary image $f(\partial\Dsmall)$
  exhibits the behavior shown in Figure \ref{fig:neglambda}(d).
\label{fig:k1graph} }
\end{figure}

\medskip
\noindent\textit{Calculation based on $t$ and $\kappa_1$.}
Now we consider how to calculate $\lambda$ corresponding to given
values of $t$ and $\kappa_1$.  We know that
\[    \kappa_1 = \kappa\, f(1) = \kappa \frac{y_2(1)}{y_1(1)},
\]
A graph of this function of $\lambda$ is illustrated in Figure
\ref{fig:k1graph}.  Since the right edge of $f(\D)$ passes through
$1$, the radius of this edge can never be less than $1/2$. Thus we
need only look for values of $\lambda$ greater than the largest value
for which $\kappa_1=2$.
 
An obvious method of calculation is to try different values of
$\lambda$, using $(t,\lambda)$ to calculate $\kappa_1$, and by
bisection or another similar method to close in on the desired value
of $\kappa_1$.  However, one can take advantage of the equivalent
formulation $a(\lambda)=0$ where
\begin{equation}  \label{eq:kappaeqn}
 a(\lambda) =  \kappa\, y_2(1) - \kappa_1 y_1(1).
\end{equation}
The coefficients of the power series $a(\lambda)$ are found in terms
of known series $\kappa$, $y_1$, and $y_2$ in $\lambda$ (note that the
first term on the right side of (\ref{eq:kappaeqn}) is a Cauchy
product).  Of course, the calculation of $y_2(1)$ involves the
iterated integrals $\X{n}$, which were not used in the 1-parameter
algorithm.  Thus the numerical calculation of $\lambda$ reduces to
finding the zero of a polynomial with least absolute value.

\begin{figure}[!b]
\pic{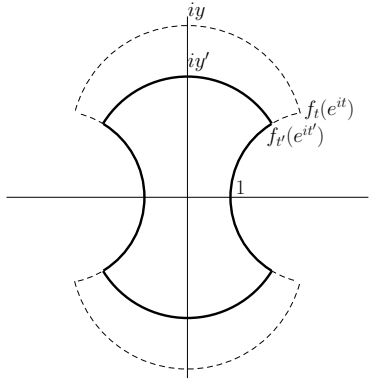}{8.5,0}{4.5}{scale=.455}
\caption{\small Increasing conformal module of $P$ corresponding
to increasing $p_2$.
 \label{fig:p2grows} }
\end{figure}

Finally we consider the upper midpoint $p_2$.

\begin{lemma} \label{lem:p2} For any fixed value of $\kappa_1\in\R$,
the geometric parameter $p_2$ is an increasing function of
$t\in(0,\pi/2)$.
\end{lemma}

\proof Fix $\kappa_2$, and let $t<t'$. The topological quadrilateral
bounded by the unit circle $\partial\D$ with vertices $\pm e^{\pm it}$
has conformal module smaller than the one with vertices $\pm e^{\pm
  it'}$. Likewise, the \scq\ $P$ whose upper edge meets the imaginary
axis in $iy$ has conformal module smaller than the one meeting in
$iy'$, when $y<y'$ (see Figure \ref{fig:p2grows}). The statement
follows.  \qed

\medskip
\noindent\textit{Two-parameter algorithm.}  We calculate
$(t,\lambda)$ from $(\kappa_1,p_2)$. Choose $t^-$ and $t^+$ near to
$0$ and $\pi/2$, respectively, so that the desired $t$ can be sought
for in the interval $t^-<t<t^+$. Begin with the midpoint
$t=(t^-+t^+)/2$, and for the three values $t^-$, $t$, $t^+$ each
paired with $\kappa_1$, find the three corresponding values of $p_2$
as described above, say $p^-$, $p$, $p^+$. By Lemma \ref{fig:p2grows},
we may repeat the process in whichever of the intervals $(p^-,p)$ or
$(p,p^+)$ the desired $p_2$ may be located.

\begin{figure}[!ht]
\pic{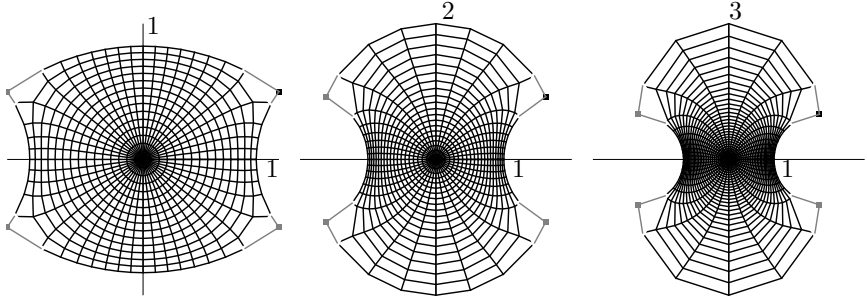}{5.5,0}{3.5}{scale=.98}
\caption{\small Images $f(\D)$ exhibiting $\kappa_1=-1$ and $p_2=1.0$,
 $2.0$, and $3.0$ respectively.
 \label{fig:2-param} }
\end{figure}

In Figure \ref{fig:2-param} some examples are given of \scq s with
prescribed values of $\kappa_1$ and $p_2$, calculated via this 
two-parameter algorithm.

\medskip
Variations on the above considerations can easily devised for finding
$(t,\lambda)$ in terms of $(\kappa_1,\kappa_2)$, etc.

\begin{figure}[!b]
\pic{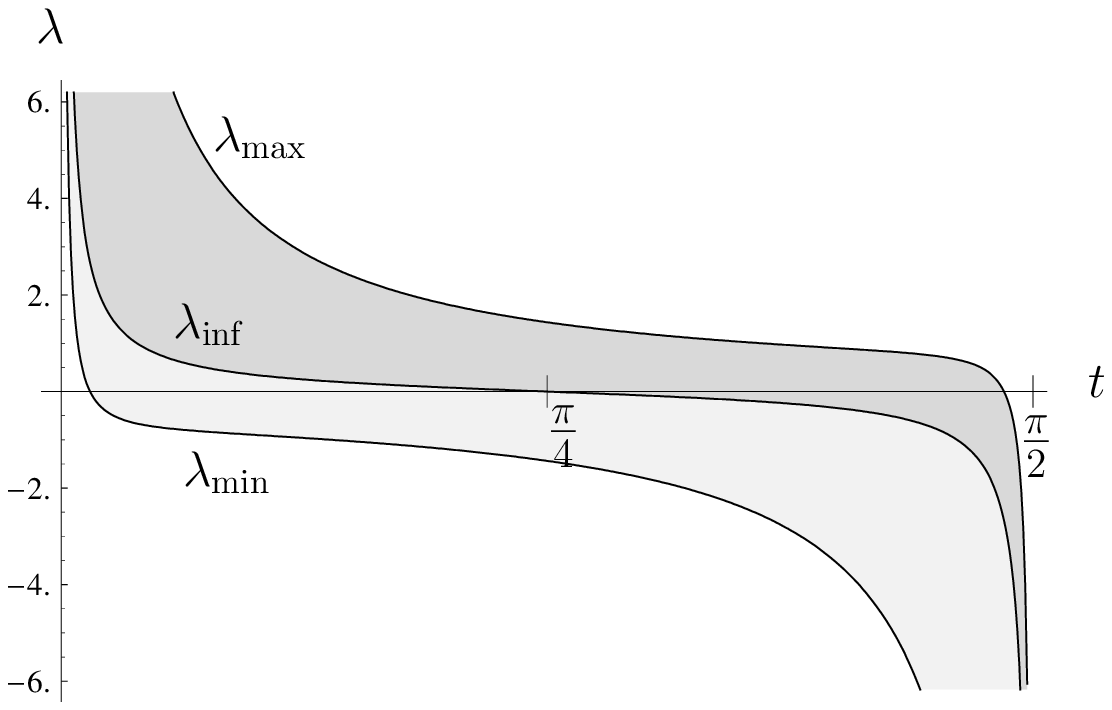}{5.5,0}{3.3}{scale=.5} 
\pic{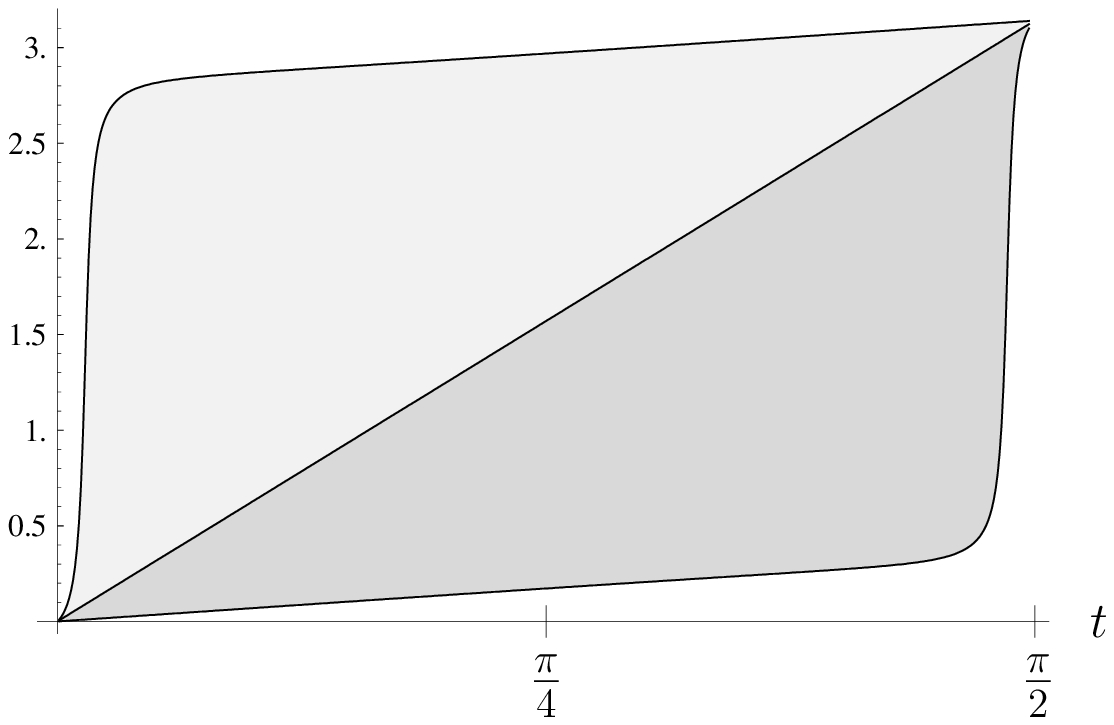}{11.,-.3}{0}{scale=.5}
\caption{\small (left) Domain of values of $(t,\lambda)$  for which
the conformal mapping to an \scq\ is univalent. (right) Domain   rescaled
by application of $\mbox{arccot}\,4\lambda$ to the vertical axis.\label{fig:univalent}  }
\end{figure}

\section{Domain of univalence.}  The above considerations
permit us to calculate easily the complete set of $(t,\lambda)$ for
which the solution to $S_f=R_{t,\lambda}$ is univalent on $\D$. From
the discussion of Figure \ref{fig:neglambda}, as $\lambda$
decreases starting from $\lambda_\infty$, the uppermost vertices
$f(e^{it})$, $f(e^{i(\pi/2-t)})$ meet on the positive imaginary axis
at a certain critical value $\lambda=\lambda_{\rm min}$.  At this
value the right and left image edges, being orthogonal to the upper
edge, must lie within the extended imaginary axis.  From this it follows
that $f(1)=\infty$, $f(-1)=-\infty$, and consequently $y_1(1)=0$.
(Note that $y_2(1)\not=0$ by linear independence.) Of course by
(\ref{eq:curvecond}) this implies $\kappa=0$, as was previously
discussed.  Similarly, as $\lambda$ increases from $\lambda_\infty$
the critical value $\lambda_{\rm max}$ approximated in Figure
\ref{fig:biglambda} is attained when $f(e^{it})$, $f(e^{-it})$ meet on
the positive real axis, and $f(i)=i\infty$, $f(-i)=-i\infty$, so
instead we solve $y_1(i)=0$.  In view of (\ref{eq:Rsymm}) this can be
solved by using $\pi/2-t$ in place of $t$.  Thus for each
$t\in(0,\pi/2)$ we have values
\[   \lambda_{\rm min}(t) < \lambda_\infty (t) <\lambda_{\rm max}(t)
\]
and in fact
\[  \lambda_{\rm min}(\pi/2-t) = -\lambda_{\rm max}(t).
\]
The domain of univalence 
$\{(t,\lambda)\colon\ \lambda_{\rm min}(t)\le t \le \lambda_{\rm max}(t)\}$ is shown in Figure \ref{fig:univalent}.

\section{Numerical Results \label{sec:numeric}}
Calculations were performed in \textit{Mathematica} 6.0 (Wolfram).

The numerical integrations were carried out by subdividing $[0,1]$
into groups of 5 consecutive intervals of length $1/M$, and
multiplying each group by a matrix equivalent to integrating the
degree-$5$ polynomial passing through the corresponding values of the
integrand.

In the calculation of the iterated integrals $I_n$ determined by
functions $q_0$, $q_1$ satisfying a bound $|q_n(z)|\le K$, from the
definition we have $|I_n(z)|\le K (\sup|I_{n-1}|)z$ for $z\in[0,1]$,
and thus $|I_n(z)|\le K^nz^n/n!$.  From this it is seen that the
Taylor coefficients of the power series we are considering tend to
zero quite rapidly.  Instead of working out here the simple a
priori error bounds which result from this, we follow the time-honored
procedure of observing how the numerical method actually works. In
terms of the two main parameters, $M$ (the number of subdivisions of
$[0,1]$ and $N$ (the degree of the polynomial approximation to the
power series), the following tables show the precision obtained for
different values of $(t,\lambda)$.  Naturally, as
$|\lambda-\lambda_\infty|$ increases, so does the numerical error
(recall that each $\lambda_\infty$ tabulated depends on the
corresponding value of $t$).

\begin{center}
Required values of $M,N$ for calculation of $\kappa$ to 5 significant
figures
\end{center} \vspace{-2ex}
\[
\begin{array}{c||r|r}
\multicolumn{3}{c}{\lambda=\lambda_\infty-1}\\
  t/\pi & M & N \\\hline
  0.05  & 40 & 15  \\
  0.08  & 25 & 15  \\
  0.1   & 15 & 15  \\
  0.2   & 10 & 15  \\
  0.3   & 10 & 15  \\
  0.4   & 10 & 11  \\
  0.45  & 10 & 11  \\
  0.48  & 10 &  9
\end{array}\quad\quad
\begin{array}{c||r|r}
\multicolumn{3}{c}{\lambda=\lambda_\infty-2}\\
  t/\pi & M & N \\\hline
  0.05  & 40 & 15  \\
  0.08  & 25 & 15  \\
  0.1   & 20 & 15  \\
  0.2   & 15 & 15  \\
  0.3   & 15 & 10  \\
  0.4   & 10 & 11  \\
  0.45  & 10 & 11  \\
  0.48  & 10 &  9
\end{array}\quad\quad
\begin{array}{c||r|r}
\multicolumn{3}{c}{\lambda=\lambda_\infty-5}\\
  t/\pi & M & N \\\hline
  0.05  & 60 & 15 \\
  0.08  & 25 & 15 \\
  0.1   & 20 & 15 \\
  0.2   & 20 & 15 \\
  0.3   & 20 & 15 \\
  0.4   & 15 & 15 \\
  0.45  & 15 & 15 \\
  0.48  & 10 & 12
\end{array}
\]

\newpage

\begin{center}
Required values of $M,N$ for calculation of $\kappa$ to 8 significant
figures
\end{center}

 \vspace{-2ex}
\[
\begin{array}{c||r|r}
\multicolumn{3}{c}{\lambda=\lambda_\infty-1}\\
  t/\pi & M & N \\\hline
  0.05  & 85 & 25  \\
  0.08  & 75 & 25  \\
  0.1   & 65 & 25  \\
  0.2   & 40 & 25   \\
  0.3   & 30 & 20  \\
  0.4   & 25 & 15  \\
  0.45  & 25 & 15  \\
  0.48  & 20 & 15 
\end{array}\quad\quad
\begin{array}{c||r|r}
\multicolumn{3}{c}{\lambda=\lambda_\infty-2}\\
  t/\pi & M & N \\\hline
  0.05  & 90 & 25  \\
  0.08  & 65 & 25  \\
  0.1   & 50 & 25  \\
  0.2   & 30 & 25  \\
  0.3   & 25 & 35  \\
  0.4   & 25 & 20  \\
  0.45  & 25 & 15  \\
  0.48  & 25 & 5 
\end{array}\quad\quad
\begin{array}{c||r|r}
\multicolumn{3}{c}{\lambda=\lambda_\infty-5}\\
  t/\pi & M & N \\\hline
  0.05  & 110 & 35  \\
  0.08  & 110 & 35  \\
  0.1   & 110 & 35  \\
  0.2   & 85 &  35 \\
  0.3   & 40 &  30 \\
  0.4   & 30 &  25 \\
  0.45  & 25 &  20 \\
  0.48  & 20 & 15
\end{array}
\]

One surprising fact is that the algorithm has no difficulty
approaching $t=\pi/2$, even though the conformal module of the
\scq\ tends to infinity and the well-known
crowding phenomenon makes graphing difficult.

In the 1-parameter algorithm, steps 1 and 2 have computational
complexity of order $O(M)$; step 3 is $O(MN)$ and step 4 is
$O(N^2)$. The solution of roots of polynomials with real coefficients
is now a well-refined science which we will not expound upon here;
one can use with goods results the internal function \texttt{NRoots}
of \textit{Mathematica}, which we believe has an operation count of
the order $O(N^2)$.  However, it is even simpler to use the standard
routine \texttt{FindRoot} for finding a single zero of an arbitrary
function, giving $\lambda_\infty$ as the starting point of the search.
Each evaluation of the function requires $O(N)$ arithmetic operations,
and the number of evaluations (by a modified form of Newton's method)
is less than proportional to the logarithm of the accuracy sought.

While the parameters $t,\lambda$ have been obtained by the SPPS
method, the calculation of the image domains in Figures
\ref{fig:biglambda}, \ref{fig:neglambda}, and \ref{fig:2-param} are
based on the function \texttt{NDSolve}, which was used to integrate
the second-order linear differential equation along radii of $\D$.
This is solely for illustration of the image domains, not for
application of the algorithms presented here.  In these figures the
vertices, being singular images, were determined in terms of 
neighboring boundary points.

The two-parameter algorithm gives no essentially new numerical
considerations.  The number of iterations of the bisection process
corresponds exactly to the desired accuracy, which also determines the
choice of $M$, $N$ as reflected in the tables above for the
calculation of $\kappa$.

\section{Conclusions}
 
 P.\ Brown \cite{Br1,Br2} studied the accessory parameters which
govern conformal mapping from the unit disk to a symmetric circular
quadrilateral, and gave algorithms for finding these parameters
in terms of geometric characteristics of the quadrilateral.

We have introduced an alternative parameter $\lambda$ in the
expression for the Schwarzian derivative of the conformal mapping, and
shown that it is identified naturally with a spectral value of a
Sturm-Liouville problem with a nonlinear boundary condition. The SPPS
method in \cite{KP} for solving Sturm-Liouville differential equations
provides simpler algorithms, especially when the solution for a
particular value of the spectral parameter is known a priori, as
is the case for the problem considered here.

One advantage of the SPPS approach as compared to other methods is
that once the simple procedure of calculating the indefinite integrals
is carried out, the problem reduces to finding roots of a polynomial.
To find the accessory parameters it is not necessary to evaluate the
conformal mapping explicitly at interior points of the domain.
Further, the method presented here does not require locating the
desired geometric characteristics approximately in a table of
previously calculated values as was done in \cite{Br2}.
  
Because of the natural relation of the second-order ordinary
differential equation to the Schwarzian derivative, we believe that
the SPPS approach given here may be applied in a similar
way to a wide variety of conformal mapping problems.

\end{document}